\documentclass[11pt]{amsart}
\usepackage{amsfonts,amssymb,amsthm,amsmath,amsxtra,amscd,verbatim,eucal}
\usepackage[all]{xy}
\usepackage[dvips]{graphics}

\theoremstyle{plain}
\newtheorem{thm}{Theorem}[section]
\newtheorem{lem}[thm]{Lemma}
\newtheorem{prop}[thm]{Proposition}
\newtheorem{cor}[thm]{Corollary}

\theoremstyle{definition}
\newtheorem{defi}[thm]{Definition}

\theoremstyle{remark}
\newtheorem{eg}[thm]{Example}

\newtheorem{rmk}[thm]{Remark}

\newtheorem*{ack}{Acknowledgements}

\numberwithin{equation}{section}

\def\N{{\mathbb N}}
\def\F{{\mathbb F}}

\def\R{{\mathbb R}}
\def\Q{{\mathbb Q}}

\def\O{\mathcal{O}}

\def\J{\mathcal{J}}
\def\I{\mathcal{I}}

\def\aa{\mathfrak{a}}
\def\bb{\mathfrak{b}}

\def\mm{\mathfrak{m}}

\def\1{{\bf 1}}
\def\ee{{\bf e}}
\def\uu{{\bf u}}
\def\vv{{\bf v}}

\def\xx{{\bf x}}
\def\zz{{\bf z}}

\def\d{\delta}

\def\p{\pi}
\def\r{\rho}
\def\s{\sigma}
\def\t{\tau}

\def\th{\theta}
\def\Th{\Theta}

\def\.{\cdot}
\def\^{\widehat}
\def\~{\widetilde}
\def\o{\circ}
\def\ov{\overline}

\def\lru{\lceil}
\def\rru{\rceil}

\def\({\left(}
\def\){\right)}

\DeclareMathOperator{\Spec} {Spec}

\DeclareMathOperator{\length} {\ell}

\DeclareMathOperator{\Vol} {Vol}
\DeclareMathOperator{\ini} {in}

\DeclareMathOperator{\lc} {lc}
\DeclareMathOperator{\Int} {Int}

\DeclareMathOperator{\Ann} {Ann}

\begin{document}

\title{Length, multiplicity, and multiplier ideals}

%\thanks{Preliminary version of June 16, 2004. Document compiled on \today.}

\author{Tommaso de Fernex}
\address{Department of Mathematics, University of Michigan,
East Hall, 525 East University Avenue, Ann Arbor, MI 48109-1109, USA}
\email{{\tt defernex@umich.edu}}

\thanks{Author's research partially supported by
the University of Michigan Rackham Research Grant and
Summer Fellowship, and by the MIUR of
the Italian Government in the framework of the National Research
Project ``Geometry on Algebraic Varieties" (Cofin 2002)}

\subjclass[2000]{Primary 14B05; Secondary 13H05, 14B07, 13H15.}
\keywords{Multiplier ideal, Samuel multiplicity, monomial ideal.}

\begin{abstract}
Let $(R,\mm)$ be an $n$-dimensional regular local ring,
essentially of finite type over a field of characteristic zero.
Given an $\mm$-primary ideal $\aa$ of $R$,
the relationship between the singularities
of the scheme defined by $\aa$ and those defined by
the multiplier ideals $\J(\aa^c)$, with $c$ varying in $\Q_+$,
are quantified in this paper by showing that
the Samuel multiplicity of $\aa$ satisfies $e(\aa) \ge (n+k)^n/c^n$
whenever $\J(\aa^c) \subseteq \mm^{k+1}$.
This formula generalizes an inequality on log canonical thresholds
previously obtained by Ein, Musta\c t\v a and the author of this paper.
A refined inequality is also shown to hold for small dimensions,
and similar results valid for a generalization of
test ideals in positive characteristics are presented.
\end{abstract}

\maketitle

\section*{Introduction}

Let $\aa$ be an ideal of an $n$-dimensional regular local ring $(R,\mm)$
essentially of finite type over a field of characteristic zero.
Consider a log resolution of $(R,\aa)$,
namely, a proper birational morphism $f$
from a smooth variety $Y$ to $\Spec R$ such that
$f^{-1}\aa \. \O_Y$ is the ideal sheaf of a simple
normal crossing divisor $E$ on $Y$.
Then, for any positive rational number $c$,
one defines the multiplier ideal
associated to the pair $(R,\aa^c)$ to be the ideal
$$
\I(\aa^c) := f_*\O_Y(K_f - [cE]) \subseteq R,
$$
where $K_f$ is the relative canonical divisor of $f$
and $[cE]$ is the integral part of the $\Q$-divisor $cE$.
This ideal is a fundamental invariant
of the singularity of the scheme defined by $\aa$, and its nice
properties related to vanishing theorems make it an important tool
in higher dimensional geometry. For
general properties of multiplier ideals, we refer to~\cite{Laz}.

For different values of $c$, the multiplier ideals $\I(\aa^c)$
offer a way to look at the singularities of the subscheme
defined by $\aa$. The goal of this paper is to
quantify this principle assuming that $\aa$ is $\mm$-primary,
by reading off conditions on the Samuel multiplicity
$e(\aa)$ and the co-length $\length(R/\aa)$ of $\aa$.

Since the Samuel multiplicity is homogeneous
with respect to integer powers of $\mm$-primary ideals, we can set
$e(\aa^c) := c^n\,e(\aa)$ for any $\mm$-primary ideal $\aa$ of
$R$ and any $c \in \Q_+$. Although this definition
is not strictly necessary, it has
the advantage of making the statement
below more symmetric and, we hope, transparent.
Using the notation just introduced and observing that $e(\mm^n) = n^n$,
the main result from~\cite{dFEM1} (namely, Theorem~1.2),
that was originally stated in terms of the log canonical threshold
of the pair $(R,\aa)$, can be reformulated as follows:
$$
e(\aa^c) \ge e(\mm^n) \quad\text{whenever}\quad \I(\aa^c) \subseteq \mm.
$$
It is natural to expect a stronger bound
on $e(\aa^c)$ when the multiplier ideal $\I(\aa^c)$
is strictly smaller than $\mm$. This motivated us towards the
next result.

\begin{thm}\label{intro-mult}
Let $R$ be as above, and let $\aa$ be an $\mm$-primary ideal $R$.
Then, for any rational $c > 0$ and any integer $k \ge 0$, we have
$$
e(\aa^c) \ge e(\mm^{n+k}) \quad\text{whenever}\quad
\I(\aa^c) \subseteq \mm^{k+1}.
$$
\end{thm}

The lower-bound on the multiplicity given in this theorem is sharp,
and the boundary cases are characterized:
equality occurs exactly when $(n+k)/c$ is an integer
and the integral closure of $\aa$ is equal to $\mm^{(n+k)/c}$
(see Theorem~\ref{equality} below).

Theorem~\ref{intro-mult} is derived by a similar result in which
the length of $R/\aa$ is bounded, instead of the multiplicity of $\aa$
(this is Theorem~\ref{length} below).
Moreover, if the dimension of $R$ is at most 3 and $c=1$,
the result on the length is, in turns, implied
by the following general formula.

\begin{thm}\label{intro-n<4}
Let $R$ be a ring as above, and assume that $R$ has dimension $n \le 3$.
Let $\aa$ be an $\mm$-primary ideal of $R$
such that $\I(\aa)$ is not trivial. Then
$$
\length(R/\aa) \ge \length\(R/\I(\aa)\mm^{n-1}\).
$$
\end{thm}

These results have analogues in positive characteristic.
In~\cite{HY}, Hara and Yoshida have introduced
a new ideal associated to a pair
$(R,\aa^c)$ when the ground field has positive characteristic.
Such ideal is denoted by $\t(\aa^c)$;
it is a generalization of the test ideal
defined by Hockster and Huneke, and can be considered as the analogue of the
multiplier ideal. In the last section of this paper we will explain
how the main theorems of this paper, stated
for multiplier ideals in characteristic zero,
also hold for this new ideal when the characteristic is positive.
The idea comes from Takagi and Watanabe's paper~\cite{TW},
where the F-pure threshold, a positive characteristic
analogue of the log canonical threshold, is defined.
Indeed it was observed in~\cite{TW} that
the same bound as the one established for log canonical thresholds
in~\cite{dFEM1}, Theorem~1.2, also holds for the F-pure threshold.

\begin{ack}
The last section of this paper was conceived
during my visit to the University of Tokyo;
I would like to thank Toshiyuki Katsura
for his invitation and support, and Shunsuke Takagi
for several enlightening discussions and remarks.
Furthermore, I would like to thank Lawrence Ein, Mel Hochster,
Rob Lazarsfeld, and Howard Thompson for many comments and
conversations that have been useful in the writing of this paper.
Finally, I would like to thank the referee for his
suggestions and corrections. 
\end{ack}

\section{Basics on monomial ideals}\label{S-monomial}

In this paper, the set of natural numbers $\N$
includes zero. For short, we denote by $\R^n_+$ the
set $(\R_{\ge 0})^n$, and view it as a subset of
the vector space $\R^n$.
Similarly, we view $\N^n$ as a subset of $\R^n_+$.
We denote by $\Vol(S)$ the volume of a measurable subset $S$ of $\R^n$.
The interior $\Int(S)$ of a subset $S \subset \R^n$ will always be
computed in $\R^n$ (even if $S$ is given as a subset of $\R^n_+$).
For a vector $\uu = (u_i) \in \R^n$,
we denote by $|\uu|$ its one-norm, that is, $|\uu| = \sum |u_i|$.
Given a hyperplane $H$ defined in $\R^n$ by an equation
of the form $\sum b_i u_i = 1$ with $b_i \in \R$,
we introduce the following notation:
$$
H^+ = \left\{ \uu \in \R^n_+ \mid \sum b_i u_i \ge 1 \right\}, \quad
H^- = \left\{ \uu \in \R^n_+ \mid \sum b_i u_i \le 1 \right\}.
$$
Let $\{\ee_1,\dots,\ee_n\}$ be the standard basis of $\R^n$,
and set $\ee = \sum \ee_i$.

Let $k$ be a field,
let $R = k[x_1,\dots,x_n]$, and denote $\mm = (x_1,\dots,x_n)$.
To a monomial $\prod x_i^{u_i} \in R$, we associate
the vector $\uu = (u_i) \in \N^n$
(we will also use the shorter notation $\xx^{\uu}$ to denote the monomial).
The Newton polygon $P(\aa)$ associated to a monomial ideal $\aa \subset R$
is, by definition, the convex hull in $\R^n_+$ of the
vectors that correspond in this fashion to the monomials in $\aa$.
Given a monomial ideal $\aa$ and a
positive integer $r$, we have $P(\aa^r) = r\,P(\aa)$.
It is natural then to set $P(\aa^c) := c\, P(\aa)$ for any positive rational
number $c$. Note that
\begin{equation}\label{c-Vol}
\Vol(\R^n_+ \setminus P(\aa^c)) = c^n\,\Vol(\R^n_+ \setminus P(\aa)).
\end{equation}
We will use the basic property saying that, for an
$\mm$-primary monomial ideal $\aa$ of $R$,
\begin{equation}\label{volume-length}
\length(R/\aa) \ge \Vol(\R^n_+ \setminus P(\aa)),
\end{equation}
and the inequality is strict if $n \ge 2$;
for a proof, see for instance~\cite{dFEM1}, Lemma~1.3.
We will also use the characterization
of Samunel multiplicity of monomial ideals that says that
\begin{equation}\label{volume-mult}
e(\aa) = n!\,\Vol(\R^n_+ \setminus P(\aa))
\end{equation}
for any $\mm$-primary monomial ideal $\aa$ of $R$ (see \cite{Tei}, page 131).

For the reminder of this section, we assume that the field $k$
has characteristic zero.

\begin{thm}[\cite{How}, Main Theorem]\label{howald}
Let $\aa$ be a monomial ideal of a polynomial ring
$k[x_1,\dots,x_n]$, where $k$ is a field of characteristic zero,
and let $c \in \Q_+$. Then $\I(\aa^c)$ is a monomial ideal as well, and
$\xx^{\uu} \in \I(\aa^c)$ if and only if
$\uu + \ee \in \Int(P(\aa^c)) \cap \N^n$.
\end{thm}

By suitably fixing a monomial order on a polynomial ring
$R = k[x_1,\dots,x_n]$,
we can define a flat deformation of an ideal $\aa$ of $R$ to its initial
ideal $\ini(\aa)$ (this is well explained in~\cite{Eis}, Section~15.8).
We can additionally assume that, for a given $c$,
this simultaneously gives a flat deformation
of $\I(\aa^c)$ to its initial ideal
$\ini(\I(\aa^c))$. Given this situation,
we will use the following property; as we will need to discuss
an analogous property in positive characteristic
(see Proposition~\ref{p-initial} below),
we include the proof for reference.

\begin{prop}[\cite{dFEM2}, Lemma~2.3]\label{initial}
Let $\aa$ be an ideal of a polynomial ring
$R = k[x_1,\dots,x_n]$, where $k$ is a field of characteristic zero,
and consider a flat deformation to initial ideals as described above. Then
$$
\I(\ini(\aa)^c) \subseteq \ini(\I(\aa^c)).
$$
\end{prop}

\begin{proof}
Let $t$ be the parameter of deformation for the degeneration to monomial
ideals fixed above.
Let $T = k[x_1,\dots,x_n,t]$, and let $\bb \subset T$ be the ideal
corresponding to the deformation of $\aa$.
There is an isomorphism from $T_t := k[x_1,\dots,x_n,t,t^{-1}]$
to $R \otimes k[t,t^{-1}]$ sending $\bb T_t$ to $\aa \otimes k[t,t^{-1}]$.
Via this isomorphism, we have
$\I((\bb T_t)^c) \cong \I(\aa^c)\otimes k[t,t^{-1}]$.
Since the family degenerating to the initial ideal is flat, we have
$$
\I(\bb^c)/((t)\cap\I(\bb^c)) \subseteq \ini(\J(\aa^c))
$$
via the canonical isomorphism $T/(t) \cong R$.
On the other hand, by standard properties of
multiplier ideals (see \cite{Laz}, Section~9.5.A), we have
$$
\I(\ini(\aa)^c) = \I((\bb/((t)\cap\bb))^c) \subseteq
\I(\bb^c)/((t)\cap\I(\bb^c)).
$$
If we put together the above inclusions, we get the assertion of the
proposition.
\end{proof}

\section{Lower-bounds to length and multiplicity}\label{S-statements}

Throughout this section, $(R,\mm)$ is an $n$-dimensional regular local ring,
essentially of finite type over a field of characteristic zero.
We begin with the following result on the co-length of
an $\mm$-primary ideal of $R$.

\begin{thm}\label{length}
Let $\aa$ be an $\mm$-primary ideal of $R$.
Fix $c \in \Q_+$ and $k \in \N$, and assume that
$\I(\aa^c) \subseteq \mm^{k+1}$. Then
\begin{equation}\label{ineq-length}
\length(R/\aa) \ge \frac{(n+k)^n}{n!\,c^n},
\end{equation}
and the inequality is strict if $n \ge 2$.
\end{thm}

We apply the above result to prove the first
theorem stated in the introduction, that we reformulate as follows.

\begin{thm}\label{mult}
With assumptions as in Theorem~\ref{length}, we have
\begin{equation}\label{ineq-mult}
e(\aa) \ge \frac{(n+k)^n}{c^n}.
\end{equation}
\end{thm}

\begin{proof}
By applying Theorem~\ref{length} to $\aa^m$ and $c/m$ in place of $\aa$ and
$c$ in the following limit, we get
$$
e(\aa) = \lim_{m \to \infty} \frac{\length(R/\aa^m)}{(m^n/n!)} \ge
\frac{(n+k)^n}{c^n}.
$$
\end{proof}

It is immediate to check that the hypotheses of Theorem~\ref{mult}
are satisfied if the integral closure of $\aa$
is a power, say $\mm^a$, of the maximal ideal $\mm$, and $c = (n+k)/a$.
Note that in this case the lower-bound to the multiplicity
given in~(\ref{ineq-mult}) is achieved. It turns out that these
are the only cases realizing the bound.

\begin{thm}\label{equality}
Under the assumptions of Theorem~\ref{mult},
$e(\aa) = (n+k)^n/c^n$ if and only if $(n+k)/c \in \N$ and
$\ov\aa = \mm^{(n+k)/c}$.
\end{thm}

The proofs of Theorems~\ref{length} and~\ref{equality}
are given in the next section. The cases $k=0$ in the three theorems above
are precisely the main results obtained in~\cite{dFEM1}
(see, respectively, Theorems~1.1,~1.2 and~1.4 ibidem).

\section{Proofs of Theorem~\ref{length}
and Theorem~\ref{equality}}\label{S-length}

The proof of Theorem~\ref{length} is conceptually elementary.
After reducing to the case of monomial ideals by a
flat degeneration, we estimate the volume
of the complement of the associated Newton polygon,
hence apply~(\ref{volume-length}). The volume is opportunely bounded
by replacing the above region by a new, smaller, symmetric
region. Roughly speaking, the new region is constructed as
the ``average" of the original one
under the action of the symmetric group on $\R^n_+$.
A key point will be that the new region we construct
is the complement (in $\R^n_+$) of a convex set.

Before passing to the proof of the theorem,
we fix some notation and prove a preliminary property on convexity.
Recall that $\R^n_+$ denotes the set $(\R_{\ge 0})^n$,
and it is viewed as a subset of $\R^n$.
The span of any non-zero vector of $\R^n_+$ is called a
{\it ray} of $\R^n_+$, and a connected bounded subset
of a ray of $\R^n_+$ containing the origin of $\R^n_+$
is said to be a {\it truncated ray}. Given any two subsets
$Q_1, Q_2 \subseteq \R^n_+$ intersecting every ray of $\R^n_+$,
we define their {\it radial sum} as the set
$$
Q_1 \star Q_2 := \bigcup_W ((Q_1 \cap W) + (Q_2 \cap W)),
$$
where the union runs over all rays $W \subset \R^n_+$ and
the sum appearing in the right hand side is the usual sum
of two subsets of a vector space. To denote a radial
summation, we will use the symbol $\bigstar$.
The following picture gives a graphical
example of radial sum of two segments $Q_1$ and $Q_2$ in $\R^2_+$.
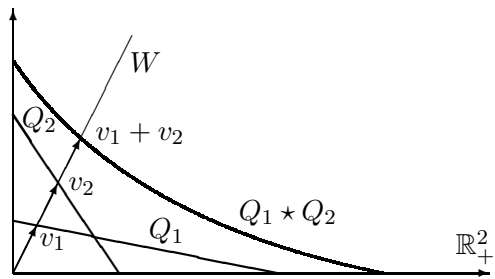
$$
\begin{picture}(180,120)
\put(0,10){\vector(1,0){180}}
\put(0,10){\vector(0,1){100}}
\put(0,10){\line(1,2){45}}
\put(0,10){\vector(1,2){9}}
\put(0,10){\vector(1,2){17}}
\put(0,10){\vector(1,2){25.5}}
\thicklines
\put(40,10){\line(-2,3){40}}
\put(100,10){\line(-5,1){100}}
\qbezier(0,90)(40,30)(140,10)
\put(15,22){\makebox(0,0){$v_1$}}
\put(26,43){\makebox(0,0){$v_2$}}
\put(48,63){\makebox(0,0){$v_1 + v_2$}}
\put(58,26){\makebox(0,0){$Q_1$}}
\put(10,68){\makebox(0,0){$Q_2$}}
\put(104,33){\makebox(0,0){$Q_1 \star Q_2$}}
\put(50,90){\makebox(0,0){$W$}}
\put(175,20){\makebox(0,0){$\R^2_+$}}
\end{picture}
$$

We will consider subsets $P \subset \R^n_+$ such that,
for any ray $W \subset \R^n_+$,
$W \setminus P$ is a truncated ray. To fix terminology, we
will say that $P$ is a {\it bounding subset} of $\R^n_+$ if
it is a subset satisfying the above property.
Newton polygons associated to $\mm$-primary monomial
ideals are examples of (closed and convex)
bounding subsets of $\R^n_+$.

\begin{prop}\label{convexity}
If $P_1,P_2 \subset \R^n_+$ are two closed, convex, bounding subsets
of $\R^n_+$,
then $P_1 \star P_2$ is also a closed, convex, bounding subset of $\R^n_+$.
\end{prop}

\begin{proof}
For any ray $W$, the set
$W \setminus (P_1 \star P_2)$ is a sum of two truncated rays
supported by $W$, hence it is a truncated ray.
In particular, $P_1 \star P_2$ is a bounding subset of $\R^n_+$.
It is also clear that $P_1 \star P_2$ is closed.
Then, to prove that $P_1 \star P_2$ is convex, it is enough to
show that it contains
the segment joining any two distinct points on
its boundary (boundaries are computed inside $\R^n_+$).
Let $\uu$ and $\vv$ be two distinct vectors lying on the
boundary of $P_1 \star P_2$. Then
$\uu = \uu_1 + \uu_2$ and $\vv = \vv_1 + \vv_2$
where, for each $i$, the vector $\uu_i$ lies on the ray spanned by $\uu$,
the vector $\vv_i$ lies on the ray spanned by $\vv$, and both
$\uu_i$ and $\vv_i$ are in the boundary of $P_i$.
Since $P_i$ is convex,
the segment joining $\uu_i$ and $\vv_i$ is contained
in $P_i$. Moreover, since $P_i$ is a bounding subset of $\R^n_+$ and
the two points are in its boundary,
the line of $\R^n$ passing through $\uu_i$ and $\vv_i$
intersects $\R^n_+$ in a bounded segment.
Note that the two segments obtained in this way
(one for each value of $i$) lie
on the two-dimensional subspace of $\R^n$ generated by
$\uu$ and $\vv$. Then, after restricting to such subspace,
the convexity of $P_1 \star P_2$ is manifestly
established by the following property:
For any two lines in $\R^2$ cutting
two non-trivial and bounded segments $L_1$ and $L_2$ on $\R^2_+$,
$L_1 \star L_2$ is a curve with concavity towards the
unbounded component of its complement in $\R^2_+$.
(This property can be easily checked by parameterizing the points on
each line by the slope of the ray they generate.)
\end{proof}

\begin{proof}[Proof of Theorem~\ref{length}]
Multipliers ideals commute with localization
and completion. Thus,
by passing to completion, hence extending zero dimensional ideals and
localizing again, we can reduce to the case in which $R$
is the localization at the origin of a polynomial
ring (for more details on the reduction,
we refer to the proof of~\cite{dFEM1}, Theorem~1.1).
After fixing a monomial order on the coordinates, we obtain
a flat degeneration of $\aa$ and $\I(\aa^c)$ to their initial ideals
$\ini(\aa)$ and $\ini(\I(\aa^c))$.
By flatness, we have $\length(R/\aa) = \length(R/\ini(\aa))$. On the other hand,
using the hypothesis that $\I(\aa^c) \subseteq \mm^{k+1}$
together with Proposition~\ref{initial}, we get
\begin{equation}\label{subsets}
\I(\ini(\aa)^c) \subseteq \ini(\I(\aa^c)) \subseteq
\ini(\mm^{k+1}) = \mm^{k+1}.
\end{equation}
Therefore, if we can prove the theorem for $\ini(\aa)$, then we can
deduce the bound for $\aa$. Thus we assume
henceforth that $\aa$ is a monomial ideal.

Let us first discuss the special case in which $k$ is a multiple of $n$, as in this
case the proof is very simple. We have $k = nq$ for some
$q \in \N$. Since $|q\ee| < k + 1$,~(\ref{subsets}) implies that
$$
\xx^{q\ee} \not \in \I(\aa^c).
$$
Then $(q+1)\ee \not \in \Int(P(\aa^c))$ by Theorem~\ref{howald}, hence
$\ee \not \in \Int(P(\aa^{c/(q+1)}))$.
Again by Theorem~\ref{howald}, this implies that
$\lc(R,\aa) \le c/(q+1)$. Therefore, applying
Theorem~1.1 from~\cite{dFEM1} to $\aa$,
we get~(\ref{ineq-length}).

Now we consider the general case. By~(\ref{volume-length})
and~(\ref{c-Vol}), it is enough to prove that
\begin{equation}\label{vol}
\Vol(\R^n_+ \setminus P(\aa^c)) \ge \frac{(n+k)^n}{n!}.
\end{equation}
In order to obtain this bound on the volume, we replace $P(\aa^c)$
by a symmetric subset of $\R^n_+$ with smaller complement.
Let $S_n$, the symmetric group of $n$ letters,
act in the obvious way on $\R^n_+$, and define
$$
Q := \frac{1}{n!}\(\mathop{\bigstar}_{\s \in S_n} \s P(\aa^c)\).
$$
Clearly, $Q$ is a $S_n$-invariant subset of $\R^n_+$.

\begin{lem}\label{Q}
$Q$ is convex, and
$\Vol(\R^n_+ \setminus Q) \le \Vol(\R^n_+ \setminus P(\aa^c))$.
\end{lem}

\begin{proof}
For short, let us denote $P(\aa^c)$ by $P$.
By recursively applying Lemma~\ref{convexity}, we see that
$\mathop{\bigstar}_{\s \in S_n} \s P$ is convex; hence
$Q$ is convex.
To prove the second part of the lemma, we fix spherical coordinates
$(\th,\r)$ on $\R^n_+$, with $\th = (\th_1,\dots,\th_{n-1})$ and
$0 \le \th_j \le \p/2$.
We denote by $\uu(\th,\r)$ the vector with spherical coordinates $(\th,\r)$.
Conversely, for a vector $\uu \in \R^n_+$, we denote by
$\th(\uu)$ and $\r(\uu)$
the multi-angle and length of $\uu$.
For any $\th \in [0,\p/2]^{n-1}$ and $\s \in S_n$, we define
$r_{\s}(\th) := \inf \{ \r \mid \uu(\th,\r) \in \s P \}$ and
$r(\th) := \inf \{ \r \mid \uu(\th,\r) \in Q \}$.
By the definition of $Q$, we have
\begin{equation}\label{rho}
r(\th) = \frac 1{n!} \sum_{\s \in S_n} r_{\s}(\th).
\end{equation}
Let $V = \{ \uu \in \R^n_+ \mid u_1 \ge \dots \ge u_n \ge 0 \}$.
This is a subcone of $\R^n_+$, and
$\R^n_+$ is the union of the cones $\s V$ as $\s$ varies in $S_n$.
Moreover, this union is almost-everywhere disjoint (the overlapping
only occurs along the boundaries of the cones $\s V$), and we clearly
have $\Vol(\s V \setminus Q) = \Vol(V \setminus Q)$ and
$\Vol(\s V \setminus P) = \Vol(V \setminus \s^{-1}P)$
for every $\s \in S_n$. Then, denoting
$\Th = \{\th \in [0,\p/2]^{n-1} \mid \uu(\th,1) \in V\}$, we have
\begin{align*}
\Vol(\R^n_+ \setminus Q) &= n!\, \Vol(V \setminus Q)
= n!\,\int_{\Th}\int_0^{r(\th)}\r^{n-1}\prod (\cos \th_i)^{i-1}\,d\r\,d\th \\
&= \int_{\Th} \(n!\,\frac{r(\th)^n}{n}\)\prod (\cos \th_i)^{i-1}\,d\th
\end{align*}
and
\begin{align*}
\Vol(\R^n_+ \setminus P) &= \sum_{\s \in S_n} \Vol(V \setminus \s P)
=\sum_{\s \in S_n} \int_{\Th} \int_0^{r_{\s}(\th)}
\r^{n-1}\prod (\cos \th_i)^{i-1}\,d\r \,d\th \\
&=\int_{\Th} \(\sum_{\s \in S_n} \frac{r_{\s}(\th)^n}{n}\)
\prod (\cos \th_i)^{i-1}\,d\th.
\end{align*}
Therefore, after we substitute the expression of $\r(\th)$
as given in~(\ref{rho}) in the integral
computing $\Vol(\R^n_+ \setminus Q)$,
the stated inequality between the two volumes
is implied by the following numerical property:
For any positive integer $n$ and positive real numbers $a_1,\dots,a_d$,
we have $d^{n-1}\sum_{j=1}^d a_j^n \ge ( \sum_{j=1}^d a_j)^n$.
A simple way to verify this property is to define
$F(\zz) = d^{n-1}\sum_{j=1}^d z_j^n - ( \sum_{j=1}^d z_j)^n$
and, denoting by $\zz_h$ the vector with $j$-component
equal to $a_j$ if $j < h$ and to $a_h$ if $j \ge h$,
use induction on $h$ to show that $F(\zz_h) \ge 0$
for every $h \in \{1,\dots,d\}$.
\end{proof}

Back to the proof of Theorem~\ref{length}, let us write
$k = nq + r$, with $0 \le r \le n-1$
(if $r=0$, then we are in the situation previously considered,
so one may assume $r >0$).
Let $J$ vary among all subsets of $\{1,\dots,n\}$
of cardinality $r$. For any such $J$, let
$\uu_J = (u_{J,i})$ be the vector whose component $u_{J,i}$ equals
$q + 2$ if $i \in J$, and $q + 1$ if $i \not \in J$.
Note that $\uu_J - \ee$ is in $\N^n$ and has one-norm
$|\uu_J - \ee| = nq + r < k+1$, hence, in view of~(\ref{subsets}),
$$
\xx^{\uu_J - \ee} \not\in \I(\aa^c).
$$
Therefore, by Theorem~\ref{howald}, $\uu_J \not \in \Int(P(\aa^c))$ for every $J$.
Then, since the set $\{ \uu_J \mid J \subset \{1,\dots,n\}, |J| = r \}$
is $S_n$-invariant, we have $\uu_J \not \in \Int(\s P(\aa^c))$
for every $J$ and every $\s \in S_n$, hence
$$
n!\, \uu_J \not \in \mathop{\bigstar}_{\s \in S_n} \Int(\s P(\aa^c))
\quad\text{for every $J$.}
$$
In conclusion, $\uu_J \not\in \Int(Q)$ for every $J$.

From now on, we fix $J = \{1,\dots,r\}$. By the convexity of $Q$, we can
find a hyperplane $H$, not containing the origin of $\R^n_+$, such that
\begin{equation}\label{H}
\uu_J \in H \quad\text{and}\quad Q \subseteq H^+.
\end{equation}
We can furthermore assume that $H$ is invariant under the action
of the stabilizer of $J$. Then $H$ has equation of the form
\begin{equation}\label{eq-H}
\sum_{i=1}^r \frac{u_i}a + \sum_{i=r+1}^n \frac{u_i}b = 1.
\end{equation}
Note that $a$ and $b$ are positive numbers, since $\R^n_+ \setminus Q$
is bounded, and we actually have
\begin{equation}\label{b}
b = \frac{a(n-r)(q+1)}{a-r(q+2)}
\end{equation}
by the first condition in~(\ref{H}).

Suppose for the moment that $a \le n+k$. It is easy to see that,
if $a = n+k$, then $H$ intersects the ray spanned by $\ee$ at the point
$((n+k)/n)\ee$. On the other hand, the value of the
left hand side of~(\ref{eq-H}) at a point $\d\ee$ is equal to
$\d(1 - r/a)/(q+1)$, and this is an increasing function of $a$.
Then, using the inequality between the arithmetic mean
and the geometric mean of the coefficients appearing
in the left hand side of~(\ref{eq-H}), we easily obtain
$$
\Vol(\R^n \cap H^-) = \frac{a^rb^{n-r}}{n!} \ge \frac{(n+k)^n}{n!}.
$$
Since $Q \subseteq H^+$, this implies~(\ref{vol})
when $a \le n+k$. Thus, we can henceforth assume that
$a \ge n+k$.

Let $V \subset \R^n_+$ be defined as in the proof of Lemma~\ref{Q}.
We remark that $\Vol(\R^n_+ \setminus Q) = n!\,\Vol(V \setminus Q)$
and $\Vol(V \setminus Q) \ge \Vol(V \cap H^-)$.
Then, in view of Lemma~\ref{Q}, in order to prove~(\ref{vol}) it suffices to
show that
$$
n!\,\Vol(V \cap H^-) \ge \frac{(n+k)^n}{n!}.
$$
The left hand side depends on $H$, that is, on $a$ and $b$, hence it
becomes a function $f(a)$ of $a$
after we substitute the expression given
in~(\ref{b}) for $b$. This function can be
computed as follows. With the help of the linear transformation
$T : \R^n \to \R^n$ given by
$u_i = \frac{w_i}i + \dots + \frac{w_n}{n}$ for $i = 1,\dots, n$,
we compute
$$
n!\,\Vol(V \cap H^-) =
\frac{1}{n!}\, a^r \prod_{i=r+1}^n \(\frac{r}{ai} + \dfrac{i-r}{bi}\)^{-1}.
$$
Substituting~(\ref{b}) for $b$, we have
$$
\(\frac{r}{ai} + \dfrac{i-r}{bi}\)^{-1} =
\frac{i(n-r)(q+1)}{i-r}\.\frac{a}{a+K_i},
\quad\text{where}\quad K_i = \frac{r(n-r)(q+1)}{i-r} - r(q+2).
$$
Therefore we obtain
$$
f(a) = C\, a^r\prod_{i=r+1}^n \frac{a}{a+K_i},
$$
where $C$ is a positive constant depending on $n$ and $k$. Deriving, we get
$$
f'(a) = \(r + \sum_{i=r+1}^n \frac{K_i}{a+K_i}\) C \, a^{r-1}
\prod_{j=r+1}^n \frac{a}{a+K_j}.
$$
It is enough to focus on the first factor in the expression of $f'(a)$.
We observe that $K_i$ get smaller as $i$ increases, and that $K_n = -r$;
hence $K_i \ge -r$ for $r+1 \le i \le n$.
Moreover, the expression $K_i/(a+K_i)$ is an increasing function
of $K_i$, for $K_i > -a$ (note that this inequality is always
satisfied). Then it is easy to see that $f'(a) > 0$ for $a > n+k$.
Since $f(n+k) = (n+k)^n/n!$, we conclude that
$f(a) > (n+k)^n/n!$ for any such $a$.
This completes the proof of the theorem.
\end{proof}

\begin{rmk}
As suggested to the author by Howard Thompson,
one can also seek a proof of Theorem~\ref{length} by using a Veronese
embedding of $R$, so to force the divisibility of $k$.
This approach requires the extension of the main
result of~\cite{dFEM1} to certain singular rings, which
is possible by the generalization of
Howald's theorem to singular Gorenstein toric rings
(for such generalization, see~\cite{Bli} or~\cite{HY}).
We prefer to work with regular rings, as this facilitates
the extension of the proof of Theorem~\ref{length} to the analogous result in
positive characteristic (see Section~\ref{S-p} below).
\end{rmk}

\begin{proof}[Proof of Theorem~\ref{equality}]
As the proof is very similar to the one of~\cite{dFEM1}, Theorem~1.4,
we will just give a brief outline, explaining the necessary modifications.
One direction is clear, so we can focus on the proof
of the ``only if'' part.

We can reduce to the case in which $R$ is the localization
at the origin of a polynomial ring.
Assuming that $e(\aa) = (n+k)^n/c^n$, we easily observe that $(n+k)/c \in \N$.
Then the theorem will follow if we show that
$\aa \subseteq \mm^{(n+k)/c}$.

We deform all the powers $\aa^r$ to monomial ideals ${\rm in}(\aa^r)$,
as follows. We first deform $\aa^r$ to the tangent cone,
and then take the initial ideal of the resulting ideal
with respect to a fixed monomial order.
In particular, doing in this way, we have $\aa \subseteq \mm^{(n+k)/c}$
if and only if $\ini(\aa) \subseteq \mm^{(n+k)/c}$.
Since we are assuming that $e(\aa) = (n+k)^n/c^n$, we have
$$
\lim_{m \to \infty} \frac{n!\, \length(R/\ini(\aa^m))}{m^n} =
\lim_{m \to \infty} \frac{n!\, \length(R/\aa^m)}{m^n} =
e(\aa) = \frac{(n+k)^n}{c^n}
$$
Set $P_m = \frac 1m\, P(\ini(\aa^m))$ and
$P_{\infty} = \bigcup_{s \in \N} P_{2^s}$. We see
from~(\ref{volume-mult}) that the volume of
$\R^n_+ \setminus P_m$ is asymptotic to $\length(R/\ini(\aa^m))/m^n$
for $m \to \infty$.
Observing that $P_{2^s} \subseteq P_{2^t}$ if $s < t$, we conclude that
$P_{\infty}$ is convex and
\begin{equation}\label{infty}
n! \Vol(\R_+^n\setminus P_{\infty}) = \frac{(n+k)^n}{c^n}.
\end{equation}
By Proposition~\ref{initial}, we have $\J(R,\ini(\aa)^c) \subseteq \mm^{k+1}$.
Then, using the inclusion $\ini(\aa)^m \subseteq \ini(\aa^m)$, we get
$\J(R,\ini(\aa^m)^{c/m}) \subseteq \mm^{k+1}$. Therefore, arguing as in the
proof of Theorem~\ref{length} and using the notation there introduced,
we see that $\uu_J \not \in \Int P(\ini(\aa^m)^{c/m})= \Int(c\,P_m)$
for every $J$. Passing to the limit
for $m = 2^s \to \infty$, we conclude that
$\uu_J \not \in \Int(c\,P_{\infty})$. Then, by the same
arguments used in the proof of Theorem~\ref{length}, applied this time
to $c\,P_{\infty}$ in place of $P(\aa^c)$,
one can check that~(\ref{infty}) can only be satisfied if
$$
\frac{n+k}{n} \. \ee \not\in \Int (c\,P_{\infty}).
$$
After this observation, we can conclude
as in the proof of~\cite{dFEM1}, Theorem~1.4.
\end{proof}

\section{A general formula on lengths in small dimension}\label{S-n<4}

We expect that the inequalities stated in Theorems~\ref{length}
and~\ref{mult} are particular cases of more
general formulae, at least when $c=1$.
We show that this is the case for the second theorem if the dimension
of $R$ is at most 3.

\begin{thm}\label{n<4}
Let $\aa$ be an $\mm$-primary ideal of a regular local ring $(R,\mm)$
of dimension $n \le 3$, essentially of finite type over a field
of characteristic zero.
Assume that $\I(\aa)$ is not trivial. Then
\begin{equation}\label{ineq-n<4}
\length(R/\aa) \ge \length(R/\I(\aa)\mm^{n-1}).
\end{equation}
\end{thm}

\begin{rmk}
Since $\I(\ov \aa) = \I(\aa)$, the statement above can be
strengthened by replacing $\length(R/\aa)$ by $\length(R/\ov \aa)$
in the left hand side of~(\ref{ineq-n<4}).
\end{rmk}

\begin{cor}\label{cor}
Keeping the notation as in Theorem~\ref{n<4},
suppose that $\I(\aa) \subseteq \mm^{k+1}$ for some $k \in \N$. Then
\begin{equation}\label{binomial}
\length(R/\aa) \ge \length(R/\mm^{n+k}) = \binom{2n+k-1}{n}.
\end{equation}
\end{cor}

\begin{rmk}
Since the binomial on the right hand side of~(\ref{binomial})
is strictly larger than $(n+k)^n/n!$ if $n \ge 2$,
Corollary~\ref{cor} gives a more precise bound on the length
than Theorem~\ref{length} when $n \in \{2,3\}$ and $c=1$.
Of course this improvement is lost in the limit computing the Samuel
multiplicity of $\aa$ (in fact we know that the statement of Theorem~\ref{mult}
is sharp).
\end{rmk}

Before passing to the proof of Theorem~\ref{n<4}, a few comments on
inequality~(\ref{ineq-n<4}). Such inequality looks very similar to
another inequality, namely
\begin{equation}\label{similar}
\length\(R/\aa\) \ge \length(R/\I(\aa\mm^{n-1})).
\end{equation}
It should be clear that~(\ref{similar}) is actually
a weaker inequality, merely following
by an inclusion of ideals. To see this, let $X = \Spec R$,
let $f : X' \to X$ be a log resolution of $(R,\aa)$, and write
$f^{-1}\mm\.\O_{X'} = \O_{X'}(-F)$ and $f^{-1}\aa\.\O_{X'} = \O_{X'}(-E)$.
Since $(R,\mm^{n-1})$ is a canonical pair,
the divisor $K_f-(n-1)F$ is effective. Thus, we have
$$
\aa \subseteq \ov{\aa} = f_*O_{X'}(-E) \subseteq
f_*O_{X'}(K_f-E - (n-1)F) = \I(\aa\mm^{n-1}).
$$
This gives~(\ref{similar}).
On the contrary, the example below shows that~(\ref{ineq-n<4})
does not follow, in general, by an inclusion of ideals.

\begin{eg}
Let $\aa$ be the ideal $(x^5,y^4,z^2)$ of $R = k[x,y,z]$
(or its integral closure). Then $\I(\aa)$ is
non-trivial, and $\aa \not \subseteq \I(\aa)\mm^2$, since
$z^2 \not \in \I(\aa)\mm^2$.
Nevertheless, it is easy to verify that
$\length(R/\aa) \ge \length(R/\I(\aa)\mm^2)$,
as we know it must be by Theorem~\ref{n<4}.
\end{eg}

The next lemma will be used to reduce the proof of Theorem~\ref{n<4}
to the case of monomial ideals.

\begin{lem}\label{reduction}
Let $R = k[x_1,\dots,x_n]_{\mm}$, with $\mm = (x_1,\dots,x_n)$.
Then, to prove a statement of the form
$$
\length(R/\bb_1) \ge \length(R/\I(\bb_2)\bb_3)
$$
where $\bb_i$ are $\mm$-primary ideals of $R$, it is enough to show that,
after computing initial ideals with respect to some monomial order on $R$,
$$
\length(R/\ini(\bb_1)) \ge \length(R/\I(\ini(\bb_2))\ini(\bb_3)).
$$
\end{lem}

\begin{proof}
We have
$$
\I(\ini(\bb_2))\ini(\bb_3) \subseteq
\ini(\I(\bb_2))\ini(\bb_3) \subseteq
\ini(\I(\bb_2)\bb_3),
$$
where the first inclusion follows from Proposition~\ref{initial},
and the second one is obvious. We get
$\length(R/\I(\ini(\bb_2))\ini(\bb_3))\ge\length(R/\ini(\I(\bb_2)\bb_3))$.
On the other hand, we have
$\length(R/\bb_1) = \length(R/\ini(\bb_1))$ and
$\length(R/\I(\bb_2)\bb_3)=\length(R/\ini(\I(\bb_2)\bb_3))$ by flatness.
\end{proof}

\begin{proof}[Proof of Theorem~\ref{n<4}]
We reduce to the case $R = k[x_1,\dots,x_n]_{\mm}$, with
$\mm = (x_1,\dots,x_n)$,
as in the proof of Theorem~\ref{length}, and
fix a monomial order inducing a flat degeneration
of the ideals $\aa$, $\I(\aa)$ and $\I(\aa)\mm^{n-1}$ to
their initial ideals $\ini(\aa)$, $\ini(\I(\aa))$ and
$\ini(\I(\aa)\mm^{n-1})$.
By the semi-continuity of the multiplier ideal, $\I(\ini(\aa)) \ne R$
if $\I(\aa) \ne R$. Therefore $\ini(\aa)$
satisfies the assumption of the theorem.
Then, by Lemma~\ref{reduction}, it is enough to
show that~(\ref{ineq-n<4}) holds for $\ini(\aa)$
in place of $\aa$. So, we can assume
that $\aa$ is a monomial ideal.

The theorem is trivial if $n=1$; we will discuss
independently the two cases when $n$ is 2 or 3.
Suppose that $n=2$.
Fix $\uu \in \N^2$ such that $\xx^{\uu} \in \aa$.
Note that $\uu \ne (0,0)$, since $\aa \ne R$.
Thus there is a choice of $i \in \{1,2\}$ for which
$\uu + \ee_i \in \Int(P(\aa))$.
This gives
$$
\xx^{\uu + \ee_i - \ee} \in \I(\aa)
$$
by Theorem~\ref{howald}. Then, taking $\{i,j\} = \{1,2\}$, we get
$$
\xx^{\uu} = \xx^{\uu + \ee_i - \ee}\.\xx^{\ee_j} \in \I(\aa)\mm.
$$
We conclude that $\aa \subseteq \I(\aa)\mm$, and this obviously
implies the inequality between the lengths.

Assume now that $n=3$, and fix $\uu = (u_1,u_2,u_3) \in \N^3$ such that
$\xx^{\uu} \in \aa$. Note that $\uu \ne (0,0,0)$.
If $\#\{i\mid u_i = 0\} \le 1$, then
essentially the same arguments used for the two-dimensional case
give $\xx^{\uu} \in \I(\aa)\mm^2$.
It remains to analyze the case when $\#\{i\mid u_i = 0\} = 2$.
For each $i \in \{1,2,3\}$, let
$$
a_i = \min\{a \in \N \mid x_i^a \in \aa\}.
$$
Since $\I(\aa) \subseteq \mm$, we have $a_i \ge 2$ for $i=1,2,3$.
We can assume that $a_1 \ge a_2 \ge a_3$.
Let $H \subset \R^3_+$ be the plane defined by $\sum u_i/a_i = 1$.
Note that $\Int(H^+) \subseteq \Int(P(\aa))$.
We consider the vectors
$$
\vv_i := (a_i-1)\ee_i + \ee_j + \ee_k, \quad \{i,j,k\} = \{1,2,3\}.
$$
Since
$$
\frac{a_1-1}{a_1} + \frac 1{a_2} + \frac 1{a_3} =
1 - \frac{1}{a_1} + \frac 1{a_2} + \frac 1{a_3} > 1,
$$
$\vv_1$ is in $\Int(H^+)$, and so in $\Int(P(\aa))$.
Then
$$
x_1^{a_1-2} = \xx^{\vv_1 - \ee} \in \I(\aa),
$$
hence $x_1^{a_1} \in \I(\aa) \mm^2$. In particular,
we have that $x_1^a \in \I(\aa) \mm^2$ for every $a \ge a_1$.
Similarly, we see that $x_2^a \in \I(\aa) \mm^2$ for every $a \ge a_2$.
For $i=3$, there are two possibilities.
If $1/a_3 < 1/a_1 + 1/a_2$, then the previous arguments
apply to this case too,
giving $x_3^a \in \I(\aa) \mm^2$ for every $a \ge a_3$,
hence  $\aa \subseteq \I(\aa)\mm^2$. However, if
\begin{equation}\label{a_i}
\frac 1{a_3} \ge \frac1{a_1} + \frac1{a_2},
\end{equation}
then we can only conclude that
$x_3^a \in \I(\aa) \mm^2$ for $a \ge a_3 + 1$.
In other words the monomial $x_3^{a_3}$, that we know
belongs to $\aa$, may be not contained in $\I(\aa)\mm^2$.
We remark that this is the only monomial in $\aa$ that does not
necessarily belong to $\I(\aa)\mm^2$.
On the other hand,~(\ref{a_i}) implies
$$
\frac{a_1-2}{a_1} + \frac 1{a_2} + \frac 1{a_3} \ge
1 - \frac{2}{a_1} + \frac 1{a_2} + \frac1{a_1} + \frac1{a_2} > 1.
$$
Note also that $a_1 \ge 2a_3 \ge 4$ by~(\ref{a_i}).
Then $(a_1 - 2)\ee_1 + \ee_2 + \ee_3 \in \Int(P(\aa))$,
hence
$$
x_1^{a_1-3} \in \I(\aa).
$$
Therefore $x_1^{a_1-1}$, which is not an element of $\aa$,
is contained in $\I(\aa)\mm^2$. In conclusion,
the number of monomials that are not contained
in $\aa$ is not smaller than the number of those that are not
in $\I(\aa)\mm^2$. This precisely means
that $\length(R/\aa) \ge \length(R/\I(\aa)\mm^2)$.
\end{proof}

\begin{rmk}
We suspect it might be possible to extend
Theorem~\ref{n<4} to all dimensions by a
suitable count of ``gains and losses''
similarly to the way done in the proof for the case $n=3$.
\end{rmk}

\section{Analogous results in positive characteristic}\label{S-p}

Let $R$ be a commutative Noetherian ring, with identity,
of characteristic $p > 0$.
We will use the letter $q$ to denote a power $p^e$ of $p$.
For any $e \in \N$, we denote by $\,^e\! R$
the ring $R$ regarded as a $R$-module via the $e$-times iterated Frobenius map
$F^e : R \to R$. For an $R$-module $M$,
we define $\F^e_R(M) := \,^e\! R \otimes_R M$, and regard it
as an $R$-module by the action of $R = \,^e\! R$ from the left.
There is an induced $e$-times iterated Frobenius map $F^e : M \to \F^e_R(M)$.
We denote by $z^q \in \F^e_R(M)$ the image by this map of an element $z \in M$.
If $N$ is an $R$-submodule of $M$, then
$N^{[q]}_M$ will denote the image
of the induced map $\F^e(N) \to \F^e(M)$.
For instance, if $I$ is an ideal of $R$,
then $I^{[q]}_R$ is the ideal of $R$
generated by the $q$-powers of the elements of $I$.

Consider an ideal $\aa$ of $R$ (we will always assume that
$\aa$ contains a nonzero divisor of $R$), and fix a positive rational number $c$.
Associated to $\aa$ and $c$, Hara and Yoshida have introduced
in~\cite{HY} a new invariant of tight closure, the definition
of which we now recall.

\begin{defi}[Hara-Yoshida]
Let $N \subseteq M$ be
$R$-modules. Given a rational number $c > 0$,
the {\it $\aa^c$-tight closure} $N^{*\aa^c}_M$ of $N$
in $M$ is defined to be the submodule of $M$ consisting
of all elements $z$ for which there exists an element $b \in R$,
not contained in any minimal prime ideal,
such that $bz^q\aa^{\lru cq \rru} \subseteq N^{[q]}_M$
for all $q = p^e \gg 0$.
\end{defi}

For the main properties of this operation, we refer to~\cite{HY}.
We limit ourselves to recall that, if $\aa = R$, then
the $R^c$-tight closure $N^{*R^c}_M$ is nothing but the usual tight closure
$N^*_M$ introduced by Hochster and Huneke in~\cite{HH}.

In~\cite{HY}, the above more general notion of tight closure
is used to define a new ideal associated to $\aa$ and $c$.
This is denoted by $\t(\aa^c)$, and is
a generalization of the test ideal introduced in~\cite{HH}.

\begin{defi}[Hara-Yoshida]
Let $E = \bigoplus_{\mm}E_R(R/\mm)$ be the direct
sum, taken over all maximal ideals $\mm$ of $R$,
of the injective envelopes of the residue
fields $R/\mm$. Then we define
$$
\t(\aa^c) := \bigcap_{M \subseteq E} \Ann_R(0^{*\aa^c}_M),
$$
where $M$ runs through all finitely generated $R$-submodules of $E$.
\end{defi}

We now collect a few properties of these ideals.
For simplicity, we henceforth assume that
$R$ is a finite dimensional regular ring essentially of finite
type over a perfect field $k$ of characteristic $p > 0$.
We refer the reader to~\cite{HY} and~\cite{HT}
for the precise statements in full generality, and for other properties
as well.

By assuming that $R$ is regular, already the definition
of the ideal $\t(\aa^c)$ becomes more
manageable, as in this case Theorem~1.13 in~\cite{HY} gives
$$
\t(\aa^c) = \Ann_R(0^{*\aa^c}_E).
$$
Hara and Yoshida also prove that, if
$\aa$ is reduced from characteristic zero
to characteristic $p \gg 0$ together with a log resolution,
then $\t(\aa^c) = \I(\aa^c)$ (see~\cite{HY}, Theorem 6.8).
More generally, the analogy between the two ideals
is expressed by several similar properties they satisfy.

It follows from the definition and the fact that
$E_R(R/\mm) = E_{R_\mm}(R_\mm/\mm R_\mm)$
that $\t((\aa R_{\mm})^c) = \t(\aa^c)R_{\mm}$
for any maximal ideal $\mm$ of $R$.
Furthermore, it is shown in~\cite{HT}, Proposition~3.2, that,
if $(R,\mm)$ is a local ring and $\^R$ is its
$\mm$-adic completion, then $\t((\aa\^R)^c) = \t(\aa^c)\^R$.
Bearing these things in mind, we can state the next properties
without assuming that $R$ is a (complete) local ring.

The ideal $\t(\aa^c)$
behaves similarly to multiplier ideals with respect to restrictions:
if $x$ is a non-zero element of $R$
and $S := R/xR$ is normal, then $\t((\aa S)^c) \subseteq \t(\aa^c)S$.
Moreover, if $R$ is a polynomial ring (or its
localization at the origin) and
$\aa$ is a monomial ideal, then the multiplier ideal
of $\aa^c$ is defined, and $\t(\aa^c) = \I(\aa^c)$;
in particular, $\t(\aa^c)$ is a monomial ideal
and, using notation as introduced in Section~\ref{S-monomial},
$\xx^{\uu} \in \t(\aa^c)$ if and only if
$\uu + \ee \in \Int(P(\aa^c)) \cap \N^n$.
For these two properties, see~\cite{HY}, Theorem~6.10.
We also recall that, if $\ov \aa$ is the integral closure
of $\aa$, then $\t(\aa^c) = \t(\ov \aa^c)$;
this follows from~\cite{HY}, Proposition~1.3.

In order to establish positive characteristic analogues
of the results proven in the previous sections of this paper,
it is necessary to control ideals of the form $\t(\aa^c)$ under degeneration
to monomial ideals.

\begin{prop}\label{p-initial}
Fix a perfect field $k$ of characteristic $p > 0$,
let $R = k[x_1,\dots,x_n]$, and consider
an ideal $\aa$ of $R$ together with a positive rational
number $c$. Then, for any monomial order
giving a flat degeneration of $\aa$ and $\t(\aa^c)$
to their initial ideals $\ini(\aa)$ and $\ini(\t(\aa^c))$,
we have
$$
\t(\ini(\aa)^c) \subseteq \ini(\t(\aa^c)).
$$
\end{prop}

The proof of Proposition~\ref{p-initial} is very similar
to that of Proposition~\ref{initial}, and we will sketch it below.
Before doing so, we only need to fix the following property.

\begin{lem}\label{pullback}
Let $R$ and $S$ be two regular rings of finite dimension, both
essentially of finite type over a perfect field $k$ of characteristic $p > 0$.
Let $T = R \otimes S$. Then, for any ideal $\aa$ of $R$,
$$
\t((\aa \otimes S)^c) = \t(\aa^c)\otimes S.
$$
\end{lem}

\begin{proof}
We can check the equality locally,
even passing to completion. Thus we may reduce to the case
when $R$ and $S$ are formal power series rings and
$T = R \^\otimes S$.
In fact, it is enough to prove the proposition
when $S = k[[t]]$, as we can then recursively apply the result
to get the statement for any dimension of $S$.

The injective envelope of the residual fields of the
rings $R$, $S$ and $T$ are respectively
$E_R = k(x_1,\dots,x_n)/k[x_1,\dots,x_n]$,
$E_S = k(t)/k[t]$, and
\begin{equation}\label{E_T}
E_T = k(x_1,\dots,x_n,t)/k[x_1,\dots,x_n,t] = E_R \otimes E_S.
\end{equation}
Moreover, we have $\F^e(E_T) = \F^e(E_R) \otimes \F^e(E_S)$.

We have $\t(((\aa \otimes S)^c)T) = \Ann_T(0^{*(\aa \otimes S)^c}_{E_T})$ and
$$
(\t(\aa^c) \otimes S)T = (\Ann_R(0^{*\aa^c}_{E_R}) \otimes E_S)T
= \Ann_T(0^{*\aa^c}_{E_R} \otimes E_S).
$$
Thus, to prove the lemma, it is enough to show that
$$
0^{*(\aa \otimes S)^c}_{E_T} = 0^{*\aa^c}_{E_R} \otimes E_S.
$$
One easily sees that the right hand side is contained in the
left hand side. To check the other inclusion, let
$z \in 0^{*(\aa \otimes S)^c}_{E_T}$.
Since $E_S = k(t)/k[t]$, we can write
$$
z = \sum_{i=1}^m w_i \otimes t^{-a_i},
$$
with $w_i \in E_R$ and $0 < a_1 < \dots < a_m$. By~\eqref{E_T}, we have
$$
w_i \otimes t^{-a_i} \in 0^{*(\aa \otimes S)^c}_{E_T}
\quad\text{for every $i$.}
$$
Observing that
$(\aa \otimes S)^{\lru cq \rru} = \aa^{\lru cq \rru}\otimes S$, this gives
$$
0 = (w_i^q \otimes t^{-a_iq})(\aa^{\lru cq \rru}\otimes S) =
(w_i^q \aa^{\lru cq \rru}) \otimes (t^{-a_iq} S)
\quad\text{in $\F^e(E_T)$}
$$
for every $i$ and all $q = p^e \gg 0$. Since
$t^{-a_iq} S$ is never zero in $\F^e(E_S)$,
we conclude that $w_i^q \aa^{\lru cq \rru} = 0$ in $\F^e(E_R)$
for sufficiently large $q$. Therefore $w_i \in 0^{*\aa^c}_{E_R}$. This in turns
yields $z \in 0^{*\aa^c}_{E_R} \otimes E_S$.
\end{proof}

\begin{proof}[Proof of Proposition~\ref{p-initial}.]
As in the proof of Proposition~\ref{initial},
let $\bb \subset T := k[x_1,\dots,x_n,t]$ be the ideal
corresponding to the deformation of $\aa$, and consider
the isomorphism from $T_t := k[x_1,\dots,x_n,t,t^{-1}]$
to $R \otimes k[t,t^{-1}]$ sending $\bb T_t$ to $\aa \otimes k[t,t^{-1}]$.
Via this isomorphism, and by Lemma~\ref{pullback}, we have
$$
\t((\bb T_t)^c) \cong \t(\aa^c)\otimes k[t,t^{-1}].
$$
From here the proof proceeds as the one of Proposition~\ref{initial}.
\end{proof}

Finally, we can state analogous results
as those proven in the previous sections for multiplier ideals.
So, let $(R,\mm)$ and be an $n$-dimensional regular local ring,
essentially of finite type over a perfect field $k$ of
positive characteristic,
and consider an $\mm$-primary ideal $\aa$ of $R$.

\begin{thm}
With the above notation,
assume that $\t(\aa^c) \subseteq \mm^{k+1}$
for some $c \in \Q_+$ and $k \in \N$. Then
$$
\length(R/\aa) \ge \frac{(n+k)^n}{n!\,c^n},
$$
with strict inequality if $n \ge 2$, and
$$
e(\aa) \ge \frac{(n+k)^n}{c^n}.
$$
Moreover, equality holds in the last formula
if and only if $(n+k)/c \in \N$ and
the integral closure of $\aa$ is equal to $\mm^{(n+k)/c}$.
\end{thm}

\begin{thm}
With the above notation, assume that $R$ has
dimension $n \le 3$ and that $\t(\aa)$ is not trivial. Then
$$
\length(R/\aa) \ge \length(R/\t(\aa)\mm^{n-1}).
$$
\end{thm}

The proof of these theorems
are the same as those of the corresponding statements
for multiplier ideals in characteristic zero.

\providecommand{\bysame}{\leavevmode \hbox \o3em
{\hrulefill}\thinspace}

\end{document}